\documentclass[11pt, leqno]{article}
\usepackage{amsmath}
\usepackage{amsfonts}
\usepackage{amssymb}
\usepackage{graphicx}
\usepackage{color}
\newtheorem{theorem}{Theorem}
\newtheorem{lemma}{Lemma}
\newtheorem{proposition}{Proposition}

\def\qed{\hfill $\square$}

\begin{document}

\title{\bf \Large A Local Curvature Bound in Ricci Flow }

\author{Peng Lu}

\date{June 04, 2009}

\maketitle

\begin{abstract}
{In this note we give a proof of a result which is closely related to Perelman's
theorem in Section 10.3 of the paper The entropy formula for the Ricci flow and
its geometric applications \cite{Pe02I}.}
\end{abstract}

\noindent \textbf{1 Introduction}. In \cite[Section 10.3]{Pe02I}  G. Perelman gives
the following theorem.
\begin{theorem} \label{Thm 10.3}
There exist $\epsilon, \, \delta >0$ with the following property.
Suppose $g_{ij}(t)$
is a smooth solution to the Ricci flow on $[0, (\epsilon r_0)^2 ]$,
and assume that at
$t=0$ we have $|\operatorname{Rm} |(x) \leq r_0^{-2}$ in $B(x_0,r_0)$, and
$\operatorname{Vol} B(x_0,r_0) \geq (1-\delta) \omega_n r_0^n$,
where $\omega_n$ is the
volume of the unit ball in $\mathbb{R}^n$.
Then the estimate $|\operatorname{Rm} |(x,t)
\leq (\epsilon r_0)^{-2}$ holds whenever $0 \leq t \leq (\epsilon r_0)^{2},
\, dist_t(x,x_0) < \epsilon r_0$.
\end{theorem}
He continues: \lq\lq The proof is a slight modification
of the proof of theorem 10.1, and is left to the reader.
A natural question is whether the assumption on the
volume of the ball is superfluous.\rq\rq

In this note by using the idea in the proof of Perelman's pseudo locality theorem
\cite[Theorem 10.1] {Pe02I}(see Theorem \ref{thm pseduo Perelman} below for the statement), we will show:
\begin{theorem} \label{thm variant}
Given $n \geq 2$ and $v_0 >0$, there exists $\epsilon_0 >0$
depending only on $n$ and $v_0$, which has the following property.
For any $r_0 >0$ and $\epsilon \in (0, \epsilon_0]$
suppose that $(M^n,g(t)), \, t \in [0, (\epsilon r_0)^2]$, is a complete smooth
solution to the Ricci flow with bounded sectional curvature, and assume that at
$t=0$ for some $x_0 \in M$
 we have curvature bound $|\operatorname{Rm} |(x,0) \leq r_0^{-2}$
 for all $x \in B_{g(0)}(x_0,r_0)$, and volume lower bound
$\operatorname{Vol}_{g(0)} \left ( B_{g(0)}(x_0,r_0) \right ) \geq v_0 r_0^n$.
Then $|\operatorname{Rm} |(x,t) \leq (\epsilon_0 r_0)^{-2}$
for all $t \in [0, (\epsilon r_0)^{2}]$ and $x \in B_{g(t)}(x_0, \epsilon_0 r_0)$.
\end{theorem}

In Sextion 2 we will give a proof of Theorem \ref{thm variant}
using two technical lemmas which will be proved in Section 3.
In Section 4 we will give two examples and a remark.
The first example shows that the curvature bound
in Theorem \ref{thm variant} is false without the assumption
$\operatorname{Vol}_{g(0)} \left ( B_{g(0)}(x_0,r_0) \right ) \geq v_0 r_0^n$.
The second example shows that the curvature bound
in Theorem \ref{thm variant} is false without the assumption that the Ricci flow
is complete.
The remark says that Theorem \ref{thm variant} follows from Theorem \ref{Thm 10.3}
and the proof of Lemma \ref{lem You Are Close}.

\vskip .1cm
\noindent \textbf{2 Proof of Theorem \ref{thm variant}}.
First we give a proof of Theorem \ref{thm variant} assuming Proposition
\ref{prop variant} below.
Then we will prove the proposition.
\begin{proposition} \label{prop variant}
Given $n \geq 2$ and $v_0 >0$, there exists $\epsilon_0 >0$ depending only
on $n$ and $v_0$
which has the following property.
For any $r_0 >0$ and $\epsilon \in (0, \epsilon_0]$
suppose that $(M^n,g(t)), \, t \in [0, \left ( \epsilon r_0 \right )^2 ]$, is a
complete smooth solution to the Ricci flow with bounded sectional curvature,
and assume that at $t=0$ for some $x_0 \in M$
we have curvature bound $|\operatorname{Rm} |(x,0) \leq r_0^{-2}$
for all $x \in B_{g(0)}(x_0,r_0)$, and volume lower bound
$\operatorname{Vol}_{g(0)} \left ( B_{g(0)}(x_0,r_0) \right ) \geq v_0 r_0^n$.
Then $|\operatorname{Rm} |(x,t)
\leq (\epsilon_0 r_0)^{-2}$ for all $t \in [0, (\epsilon r_0)^{2}]$ and
$x \in B_{g(0)} \left (x_0, e^{n-1} \epsilon_0 r_0 \right )$.
\end{proposition}

\noindent \textbf{Proof of Theorem \ref{thm variant}}.
It suffices to prove the following statement.
For the solution $g(t)$ in Proposition \ref{prop variant} we have
\begin{equation}
B_{g(t)} \left (x_0, \epsilon_0 r_0 \right ) \subset B_{g(0)}
\left (x_0, e^{n-1} \epsilon_0 r_0 \right )
\text{~ for any } t \in [0, (\epsilon r_0)^{2}]. \label{eq ball relation}
\end{equation}
We will prove (\ref{eq ball relation}) by contradiction.

If (\ref{eq ball relation}) is not true,
there is a point $x \in B_{g(t)} \left (x_0, \epsilon_0 r_0 \right )
\setminus B_{g(0)} \left (x_0, e^{n-1} \epsilon_0 r_0 \right )$.
Let $\gamma(s)$, $0 \leq s \leq s_0$, be a unit-speed minimal geodesic with
respect to metric $g(t)$ such that $\gamma(0)=x_0$ and $\gamma(s_0)=x$.
Then $s_0 < \epsilon_0 r_0$,
and there is a $s_1 \in (0,s_0]$ such that $\gamma(s_1) \in
\partial \left ( B_{g(0)} \left (x_0, e^{n-1} \epsilon_0 r_0 \right ) \right )$
and $\gamma([0,s_1)) \subset  B_{g(0)}\left (x_0, e^{n-1} \epsilon_0 r_0 \right )$.
In particular, the length satisfies
\begin{equation}
L_{g(0)} \left (\left . \gamma \right \vert_{[0,s_1]}
\right ) \geq e^{n-1} \epsilon_0 r_0. \label{eq length lower bound}
\end{equation}

From the curvature bound $|\operatorname{Rm} |(x,t) \leq ( \epsilon_0
r_0)^{-2}$ in Proposition \ref{prop variant} and the Ricci flow equation, we have
\[
\left | \gamma ^\prime (s) \right |_{g(0)} \leq e^{(n-1) }
\left  | \gamma ^\prime (s) \right |_{g(t)} \text{~ for all }
t \in [0, ( \epsilon r_0)^{2}] \text{ and } s \in [0,s_1].
\]
Hence 
\[
L_{g(0)} \left (\left . \gamma \right \vert_{[0,s_1]} \right ) \leq
\int_0^{s_1} e^{n-1} \left  | \gamma ^\prime (s) \right |_{g(t)} ds
 \leq e^{n-1 } \cdot s_0 < e^{n-1} \epsilon_0 r_0.
\]
This contradicts with (\ref{eq length lower bound}). Hence (\ref{eq ball relation})
is proved, and Theorem \ref{thm variant} is proved assuming Proposition 
\ref{prop variant}.
\qed

\vskip .1cm
In the rest of this section we give a proof of Proposition \ref{prop variant}.

\noindent \textbf{Proof of Proposition \ref{prop variant}}.
Let $\tilde{g}(t) \doteqdot (r_0)^{-2} g((r_0)^{2}t)$ be the
parabolically scaling of $g(t)$. 
The proposition holds for $g(t)$ and $r_0$ if and only if the proposition holds
 for $\tilde{g}(t)$ and $r_0 =1$. 
Hence it suffices to prove the proposition for $r_0 =1$ 
which we assume from now on. We will prove the proposition for $r_0=1$ 
by contradiction argument.

Suppose the proposition is not true.
Then there are $n \geq 2, \, v_0 >0$, a sequence of $\epsilon_{0i} \rightarrow 0^+$,
a sequence of $\epsilon_{i} \in (0, \epsilon_{0i}]$,
a sequence of complete smooth solutions to the Ricci flow
$\left (M^n_i,g_{i}(t) \right ), \, t \in [0, \epsilon^2_i]$,
with bounded sectional curvature, and a sequence of points $x_{0i} \in M_i$, 
such that the following is true for each $i$:

\noindent (i) $|\operatorname{Rm}_{g_i} |(x,0) \leq 1$ for all $x \in
B_{g_i(0)}(x_{0i},1)$. 

\noindent (ii) $\operatorname{Vol}_{g_i(0)} \left ( B_{g_i(0)}(x_{0i},1)
\right ) \geq v_0$.

\noindent (iii) There are $t_i \in (0, \epsilon_i^{2}]$ and
$x_i \in B_{g_i (0)}(x_{0i}, e^{n-1} \epsilon_{0i})$ such that
 $|\operatorname{Rm}_{g_i} |(x_i,t_i)  > \epsilon_{0i} ^{-2}$.

\noindent (iv) $\epsilon_{0i} \leq \frac{1}{8 e^{n-1}}$.

To get a contradiction from the existence of sequence $\left \{\left (M_i, g_{i}(t)
\right ) \right \}$, we need the following point-picking statement whose
proof is simpler than the proof of the point-picking claim used by Perelman
in \cite[Section 10.1]{Pe02I}.
Let $A_i \doteqdot \frac{1}{100n \epsilon_{0i}}$.

\textbf{Claim A}. Fix any $i$, there are points $(\bar{x}_i,\bar{t}_i)
\in B_{g_i(0 )}(x_{0i}, (2A_i+e^{n-1}) \epsilon_{0i}) \times(0,\epsilon_i^2]$
with $\bar{Q}_i \doteqdot | \operatorname{Rm}_{g_i} | (\bar{x}_i,
\bar{t}_i) > \epsilon_{0i}^{-2}$ such that
\[
| \operatorname{Rm}_{g_i} | (x,t) \leq 4  \bar{Q}_i
\text{ ~ for all } (x,t) \in  B_{g_i(0)} \left (\bar{x}_{i}, A_i \bar{Q}_i ^{-1/2}
\right ) \times(0,\bar{t}_i] .
\]

\textbf{Proof of Claim A}. Let $Q_i^0 \doteqdot |\operatorname{Rm}_{g_i} | ({x}_i,
{t}_i)$. If $(x_i,t_i)$ from (iii) satisfies the curvature
bound of the claim,  ie.
\[
| \operatorname{Rm}_{g_i} | (x,t) \leq 4 Q_i^0
\text{ ~for } (x,t) \in  B_{g_i(0)} \left (x_{i}, A_i (Q_i^0)^{-1/2} \right )
\times (0, t_i] ,
\]
we choose $(\bar{x}_i,\bar{t}_i) =(x_i,t_i)$ and the claim is proved.

If $(x_i,t_i)$ does not satisfy the curvature bound of the claim,
 then there is a point
\[
(x_i^1,t_i^1) \in  B_{g_i(0)} \left ({x}_{i},
 {A_i} \left ( Q_i^0 \right )^{-1/2} \right ) \times(0,t_i]
\]
such that $ | \operatorname{Rm}_{g_i} | (x_i^1, t_i^1) > 4  Q_i^0$.
We compute using $ Q_i^0 > \epsilon_{0i}^{-2}$
\begin{align*}
& d_{g_i(0)}(x_i^1,x_{0i}) \leq d_{g_i(0)}(x_i,x_{0i})
+ A_i \left ( Q_i^0 \right )^{-1/2} \\
& \leq e^{n-1} \epsilon_{0i}+ A_i \epsilon_{0i} \\
& < (2A_i +e^{n-1}) \epsilon_{0i}.
\end{align*}
If $(x_i^1,t_i^1)$ satisfies the curvature bound of the claim,
we choose $(\bar{x}_i,\bar{t}_i) =(x_i^1,t_i^1)$ and the claim is proved.

If $(x_i^1,t_i^1)$ does not satisfy the claim, let $Q_i^1 \doteqdot
 |\operatorname{Rm}_{g_i} | (x_i^1, t_i^1)$, then there is a point
\[
(x_i^2,t_i^2) \in  B_{g_i(0 )} \left (x_i^1, A_i \left( Q_i^1 \right )^{-1/2}
\right ) \times(0,t_i^1]
\]
such that $ | \operatorname{Rm}_{g_i} | (x_i^2, t_i^2) > 4  Q_i^1$.
We compute using $ Q_i^1 > 4Q_i^0$
\begin{align*}
& d_{g_i(0 )}(x_i^2,x_{0i})  \leq
d_{g_i(0 )}({x}_i^1,x_{0i}) + A_i \left ( Q_i^1\right )^{-1/2} \\
& \leq (e^{n-1}+A_i) \epsilon_{0i}+ A_i \cdot \frac{1}{2} \epsilon_{0i} \\
& <(2A_i+e^{n-1}) \epsilon_{0i}.
\end{align*}
If $(x_i^2,t_i^2)$ satisfies the curvature bound of the claim,
we choose $(\bar{x}_i,\bar{t}_i) =(x_i^2,t_i^2)$ and the claim is proved.

If $(x_i^2,t_i^2)$ does not satisfy the claim, then there will be a point
$(x_i^3,t_i^3)$ and we can continue the above process of arguments.
Hence for each $i$ either we get a finite sequence points $\{ (x_i^k,t_i^k)
\}_{k=0}^{k_i}$ where $(x_i^0,t_i^0) \doteqdot (x_i,t_i)$ such that
the claim holds by taking $(\bar{x}_i,\bar{t}_i) =(x_i^{k_i},t_i^{k_i})$,
or there is an infinite sequence of points $\{ (x_i^k,t_i^k)
\}_{k=0}^{\infty}$ which satisfies the following.
Let $Q_i^k \doteqdot |\operatorname{Rm}_{g_i} | (x_i^k, t_i^k)$,
then for each integer $k \geq 0$
\[
(x_i^{k+1},t_i^{k+1}) \in  B_{g_i(0 )} \left (x_i^k, A_i \left( Q_i^k \right )^{-1/2}
\right ) \times(0,t_i^k]
\]
such that $ | \operatorname{Rm}_{g_i} | (x_i^{k+1}, t_i^{k+1}) > 4  Q_i^k$.

Now we show that for any $i$ there can not be infinite sequence $\{ (x_i^k,t_i^k)
\}_{k=0}^{\infty}$ from which the claim follows.
We compute
\begin{align*}
& d_{g_i(0)}(x_i^{k+1},x_{0i}) \\
&\leq
d_{g_i(0 )}(x_{0i}, x_i^0) +d_{g_i(0 )}({x}_i^0,x_i^1)++d_{g_i(0 )}({x}_i^1,x_i^2)
+ \cdots + d_{g_i(0 )}({x}_i^k,x_i^{k+1}) \\
& \leq e^{n-1}\epsilon_{0i} + A_i \left( Q_i^0 \right )^{-1/2}+ A_i \left( Q_i^1
\right )^{-1/2} + \cdots +  A_i \left( Q_i^k \right)^{-1/2} \\
& \leq e^{n-1}\epsilon_{0i} + A_i \epsilon_{0i}+ A_i \frac{1}{2} \epsilon_{0i}
+\cdots + A_i \frac{1}{2^k} \epsilon_{0i} \\
& < (2A_i+e^{n-1}) \epsilon_{0i},
\end{align*}
where we have used $ Q_i^{k+1} > 4Q_i^k >4^{k+1} Q_i^0> 4^{k+1}\epsilon_{0i}^{-2}$.
For any fixed $i$, from $A_i =\frac{1}{100 n \epsilon_{0i}}$ and
$\epsilon_{0i} \leq \frac{1}{8 e^{n-1}}$,
we conclude that $(x_i^k,t_i^k)$ is in the compact set
$\overline{B_{g_i(0 )}(x_{0i}, 1)}
\times[0,\epsilon_i^2]$ for all $k$. On the other hand we have
\[
\lim_{k \rightarrow \infty }| \operatorname{Rm}_{g_i} | (x_i^{k}, t_i^{k})
\geq \lim_{k \rightarrow \infty }  4^k\epsilon_{0i}^{-2} =\infty,
\]
which is impossible. Now Claim A is proved.
\qed

\vskip .1cm
Let $(\bar{x}_i,\bar{t}_i)$ be the point given by Claim A.
We divide the rest proof of Proposition \ref{prop variant} into three cases
according to the value of
\begin{equation}
\overline{\lim}_{i \rightarrow \infty} \, \bar{t}_i \cdot
| \operatorname{Rm}_{g_i} | (\bar{x}_i, \bar{t}_i) \doteqdot \tilde{\alpha}
\label{eq tilde alpha}
\end{equation}
equals to infinite, positive finite number, or zero.
We will derive contradictions in all three cases.

\vskip .1cm
\textbf{Case 1} $\tilde{\alpha} = + \infty$.
From Claim A and the choice of $A_i =\frac{1}{100n \epsilon_{0i}}$,
by switching to a subsequence (still indexed by $i$)
we have 

\noindent (1i) $\bar{t}_i \leq \epsilon_i^2$,

\noindent  (1ii) ${\lim}_{i \rightarrow \infty} \, \bar{t}_i \cdot
| \operatorname{Rm}_{g_i} | (\bar{x}_i, \bar{t}_i) = \infty$, 

\noindent  (1iii) $d_{g_i(0)}(\bar{x}_i,x_{0i}) <\frac{1}{4}$.
In particular,
$B_{g_i(0)}(\bar{x}_i,\frac{3}{4}) \subset B_{g_i(0)}(x_{0i},1)$.

From the assumptions (i) and (ii) given at the beginning of the proof
 of Proposition \ref{prop variant} and the
Bishop-Gromov volume comparison theorem there is a constant $v_1>0$,
depending only on $n$ and $v_0$, such that
\[
\operatorname{Vol}_{g_i(0)} \left ( B_{g_i(0)}\left(x_{0i},\frac{1}{4} \right )
\right ) \geq v_1.
\]
Since the ball $B_{g_i(0)} \left(\bar{x}_i,\frac{1}{2} \right)$ contains ball
$B_{g_i(0)} \left (x_{0i},\frac{1}{4} \right)$ we have
\begin{equation}
\operatorname{Vol}_{g_i(0)} \left ( B_{g_i(0)} \left(\bar{x}_i,\frac{1}{2} \right)
\right ) \geq v_1.
\end{equation}

We define a {\em regular domain} in a smooth manifold to be a bounded domain with a
$C^1$-boundary. Recall Perelman's pseudolocality theorem \cite[Theorem 10.1]{Pe02I}
says the following (for an expository account, see, for example, Chow et al
\cite[Chapter 21]{CCGG}).
\begin{theorem}\label{thm pseduo Perelman} (Perelman) For every $\alpha > 0$
and $n \geq 2$ there exist
$\delta > 0$ and $\epsilon_0 > 0$ depending only on $\alpha$ and $n$ with the
following property. Let $(M^n, g(t)), \, t \in [0,(\epsilon r_0)^2]$,
where $\epsilon \in (0, \epsilon_0]$ and $r_0\in ( 0, \infty)$, be a complete
solution of the Ricci flow with bounded curvature and let $x_0 \in M$ be a point such that
\[
R(x, 0) \geq -r_0^{-2} \hskip .5cm \text{ for } x \in B_{g(0)}(x_0, r_0)
\]
and
\[
\left ( \operatorname{Area}_{g(0)} (\partial \Omega ) \right )^n \geq
(1- \delta)c_n \left ( \operatorname{Vol}_{g(0)}(\Omega) \right )^{n-1}
\]
for any regular domain $\Omega \subset B_{g(0)}(x_0, r_0)$, where $c_n = n^n
\omega_n$ is the Euclidean isoperimetric constant. Then we have the curvature estimate
\[
\left \vert \operatorname{Rm} \right \vert (x, t) \leq
\frac{\alpha}{t} +
\frac{1}{(\epsilon_0 r_0 )^2}
\]
for $x\in  B_{g(t)} (x_0, \epsilon_0 r_0)$ and $t \in (0, (\epsilon_0 r_0)^2]$.
\end{theorem}

Let $\delta \doteqdot \delta_0 >0$ be the constant in Theorem 
\ref{thm pseduo Perelman}
corresponding to $\alpha =1$. Applying Lemma \ref{lem You Are Close} below to
metric $4 g_i(0)$ and ball $B_{4g_i(0)}(\bar{x}_i,1) =B_{g_i(0)}\left( \bar{x}_i,
\frac{1}{2} \right )$ we conclude that there is a $r_1 <\frac{1}{2}$,
depending only on $n$, $\delta_0$ and $v_1$ but not depending on $i$, such that
\begin{equation}
\left ( \operatorname{Area}_{g_i(0)} \left(  \partial\Omega\right) \right )  ^{n}\geq\left(  1-\delta
_{0}\right)  c_{n}\left(  \operatorname{Vol}_{g_i(0)} \left(  \Omega\right)  \right)
^{n-1} \label{eq isoper appl 1}
\end{equation}
for any regular domain $\Omega\subset B_{g_i(0)}\left( \bar{x}_{i},r_1\right)$.

Let $r_2 \doteqdot \min \left \{r_1,\frac{1}{\sqrt{n(n-1)}} \right \}$,
and let $\hat{g}_i(t) =(r_2)^{-2} g_i((r_2)^2 t)$, $0 \leq t \leq (r_2)^{-2}
\epsilon_i^2$.
It follows from assumption (i) that the scalar curvature $R_{\hat{g}_i}(\cdot,0)
\geq -1$ on $B_{\hat{g}_i(0)}(\bar{x}_i,1)$.
From (\ref{eq isoper appl 1}) we have
\[
\left ( \operatorname{Area}_{\hat{g}_i(0)} \left(  \partial\Omega\right) 
\right ) ^{n}
\geq \left(  1-\delta _{0}\right)  c_{n}\left(  \operatorname{Vol}_{\hat{g}_i(0)}
\left(  \Omega\right)  \right) ^{n-1}
\]
for any regular domain $\Omega\subset B_{\hat{g}_i(0)}\left( \bar{x}_{i},1 \right)$.

For $i$ large enough we can apply Theorem \ref{thm pseduo Perelman} (using
$\alpha =1$) to $\left(
B_{\hat{g}_i(0)}\left( \bar{x}_{i},1 \right), \hat{g}_i(t) \right )$,
$0 \leq t \leq (r_2)^{-2} \epsilon_i^2$, and conclude
\[
\left | \operatorname{Rm}_{\hat{g}_i} \right |(x,t) \leq \frac{1}{t}
+\frac{1}{(r_2)^{-2} \epsilon_i^2}
\]
for $t \in (0,(r_2)^{-2} \epsilon_i^2]$ and $x \in B_{\hat{g}_i(t)}
(\bar{x}_i, (r_2)^{-1} \epsilon_i)$. Equivalently we have
\[
\left | \operatorname{Rm}_{g_i} \right |(x,t) \leq \frac{1}{t}
+\frac{1}{ \epsilon_i^2}
\]
for $t \in (0, \epsilon_i^2]$ and $x \in B_{g_i(t)}
(\bar{x}_i, \epsilon_i)$. In particular
\[
\left | \operatorname{Rm}_{g_i} \right |(\bar{x}_i,\bar{t}_i) \leq \frac{1}{\bar{t}_i}
+\frac{1}{ \epsilon_i^2} \leq \frac{2}{\bar{t}_i}
\]
for $i$ large enough. This contradicts with the assumption of Case 1 that $\tilde{\alpha}$
in (\ref{eq tilde alpha}) is infinity.

\vskip .1cm
\textbf{Case 2} $\tilde{\alpha} \in (0, \infty)$. Let $\hat{t}_i \doteqdot 
\bar{Q}_i \bar{t}_i$.
Let $\hat{g}_i(t) \doteqdot \bar{Q}_ig_i(\left ( \bar{Q}_i\right )^{-1}t)$,
$t \in [0,\hat{t}_i]$.
Let $b_0$ be a constant bigger than $\frac{11}{3}(n-1)(\tilde{\alpha}+1)+1$
to be chosen later (see (\ref{eq tem use 15}) below).
By passing to a subsequence we have

\noindent (2i) $\left\vert \operatorname*{Rm}_{\hat{g}_{i}}\right
\vert \left( x,t\right) \leq 4$ for $x\in B_{\hat{g}_{i}\left( 0\right)
}\left( \bar{x}_{i},A_i\right)$ and $t \in [0,\hat{t}_i]$,

\noindent (2ii) $\left \vert \operatorname*{Rm}_{\hat{g}_{i}} \right \vert
\left(\bar{x}_i, \hat{t}_i \right) =1 $,

\noindent (2iii) $\left\vert \operatorname*{Rm}_{\hat{g}_{i}}\right\vert
\left( x,0\right) \leq \bar{Q}_i^{-1}$
for $x\in B_{\hat{g}_{i}\left( 0\right) }\left( \bar{x}_{i},A_i\right)$,

\noindent (2iv) $\hat{t}_i \leq \tilde{\alpha}+1$, $\hat{t}_i \rightarrow \tilde{\alpha}$,
$A_i >2 e^{4(n-1)(\tilde{\alpha}+1)}b_0$, and $A_i \rightarrow \infty$.

Applying Lemma \ref{lem you choose 2} to $\hat{g}_i(t)$ with $b=b_0$ we get a function
$h_i: M_i \times [0,\hat{t}_i] \rightarrow [0,1]$ such that the support
\[
\operatorname{supp}h_i(\cdot,t) \subset \overline{B}_{\hat{g}_i(t)} \left
(\bar{x}_i, 2b_0 -\frac{11}{3}(n-1)t \right ) \subset B_{\hat{g}_{i}\left( 0\right)
}\left( \bar{x}_{i},A_i\right)
\]
and
\begin{equation*}
\left( \frac{\partial }{\partial t}-\Delta _{\hat{g}_i\left( t\right) }\right) h_i
\leq \frac{10}{b_0^2}h_i.
\end{equation*}

Recall the curvature $\operatorname{Rm}_{\hat{g}_i}$ of Ricci flow $\hat{g}_i(t)$
satisfies
\begin{equation*}
\left( \frac{\partial }{\partial t}-\Delta _{\hat{g}_i }\right)
\left\vert \operatorname{Rm}_{\hat{g}_i}\right\vert ^{2}\leq
-2\left\vert \nabla_{\hat{g}_i } \operatorname{Rm}_{\hat{g}_i}
\right\vert ^{2}+16\left\vert \operatorname{Rm}_{\hat{g}_i}\right\vert ^{3}.
\end{equation*}

Now we compute the evolution equation of $h_i\left\vert \operatorname{Rm}
_{\hat{g}_i }\right\vert ^{2}$.
\begin{align*}
&\left( \frac{\partial }{\partial t}-\Delta _{\hat{g}_i }\right)
\left( h_i\left\vert \operatorname{Rm}_{\hat{g}_i }\right\vert ^{2}\right)  \\
&=\left( \left( \frac{\partial }{\partial t}-\Delta _{\hat{g}_i
}\right) h_i \right) \left\vert \operatorname{Rm}_{\hat{g}_i } \right\vert ^{2}
+h_i\left( \left( \frac{
\partial }{\partial t}-\Delta _{\hat{g}_i }\right) \left\vert \operatorname{Rm
}_{\hat{g}_i }\right\vert ^{2}\right) -2\nabla_{\hat{g}_i } h_i\cdot
\nabla_{\hat{g}_i } \left\vert \operatorname{Rm}_{\hat{g}_i }
\right\vert ^{2} \\
&\leq \frac{10}{b_0^2}\, h_i\left\vert \operatorname{Rm}_{\hat{g}_i
}\right\vert ^{2}+h_i\left(
-2\left\vert \nabla_{\hat{g}_i } \operatorname{Rm}_{\hat{g}_i }\right\vert^{2}
+16\left\vert \operatorname{Rm}_{\hat{g}_i }
\right\vert ^{3}\right) + \frac{4 \sqrt{10}}{b_0}
\left\vert \operatorname{Rm}_{\hat{g}_i }\right\vert \cdot h_i^{1/2}
\left\vert \nabla_{\hat{g}_i } \operatorname{Rm}_{\hat{g}_i }\right\vert  \\
&\leq \left ( \frac{10}{b_0^2} + 64 \right ) \, h_i\left\vert \operatorname{Rm}_{
\hat{g}_i } \right\vert^{2}-2h_i\left\vert \nabla_{\hat{g}_i } \operatorname{Rm}
_{\hat{g}_i } \right\vert ^{2}+ \frac{16 \sqrt{10}}{b_0}
 \cdot h_i^{1/2}\left\vert
\nabla_{\hat{g}_i } \operatorname{Rm}_{\hat{g}_i }
\right\vert  \\
& \leq \left ( \frac{10}{b_0^2} + 64 \right ) \, \left( h_i\left\vert
\operatorname{Rm}_{\hat{g}_i } \right\vert ^{2}\right) +  \frac{320}{b_0^2},
\end{align*}
where we have used
\begin{equation*}
\left\vert \nabla_{\hat{g}_i } h_i \right\vert =\frac{\left\vert \phi ^{\prime }\left(
w\right) \right\vert }{b_0}\left\vert \nabla _{\hat{g}_i }d_{\hat{g}_i\left(
t\right) }\left( x,\bar{x}_i \right) \right\vert _{\hat{g}_i}
\leq \frac{ \sqrt{10}}{b_0}h_i^{1/2}
\end{equation*}
and $\left\vert \operatorname{Rm}_{\hat{g}_i }\right\vert \leq 4$ on
$\operatorname{supp}h_i(\cdot,t)$. Here $\phi$ is the function defined in the proof
of Lemma \ref{lem you choose 2}.

Let $u_i \doteqdot h_i\left\vert \operatorname{Rm}_{\hat{g}_i }\right\vert ^{2}$.
We have proved
\[
\left( \frac{\partial }{\partial t}-\Delta _{\hat{g}_i  }\right) u_i
\leq \left ( \frac{10}{b_0^2} + 64 \right ) u_i+ \frac{320}{b_0^2}
\]
on $M_i \times [0,\hat{t}_i]$.

Let $H_i>0$ be the backward heat
kernel to the conjugate heat equation on $(M_i, \hat{g}_i(t)), \,
t \in [0,\hat{t}_i]$, centered at $\bar{x}_i$, ie,
\begin{eqnarray*}
\left( \frac{\partial }{\partial t}+\Delta _{\hat{g}_i }-R_{\hat{g}_i }\right)
H_i &=&0 \\
\lim_{t\rightarrow \hat{t}_i}H_i\left( x,t\right)  &=&\delta _{\bar{x}_i}.
\end{eqnarray*}
Note that $\int_{M_i}H_i(\cdot, t) d\mu _{\hat{g}_i\left( t\right) }=1$.

Now we compute
\begin{align*}
&\frac{d}{dt}\int_{M_i}u_iH_id\mu _{\hat{g}_i } \\
&= \int_{M_i}\left( \left( \frac{\partial }{\partial t}-\Delta _{\hat{g}_i}
\right) u_i \right) H_id\mu _{\hat{g}_i }+\int_{M_i}u_i\left( \left(
\frac{\partial }{\partial t}+\Delta _{\hat{g}_i }-R_{\hat{g}_i
}\right) H_i \right) d\mu _{\hat{g}_i } \\
&\leq \int_{M_i}\left( \left ( \frac{10}{b_0^2} + 64 \right ) u_i+
\frac{320}{b_0^2} \right) H_id\mu _{\hat{g}_i } \\
&= \left ( \frac{10}{b_0^2} + 64 \right ) \int_{M_i}u_i H_i d\mu _{\hat{g}_i }
+ \frac{320}{b_0^2}.
\end{align*}
Hence it follows from a simple integration that $U_i \left( t\right) \doteqdot
\int_{M_i}u_i H_id\mu _{\hat{g}_i}$ satisfies
\begin{equation}
U_i \left( t\right) \leq e^{\left ( \frac{10}{b_0^2} + 64 \right )t}U_i
\left( 0\right) +\frac{320}{\left ( \frac{10}{b_0^2} + 64 \right )b_0^2}\left(
e^{\left ( \frac{10}{b_0^2} + 64 \right )t}-1\right) \label{eq tem use 11}
\end{equation}
for $t \in [0,\hat{t}_i]$.

By the definition of $h_i$ we have at $t=\hat{t}_i$
\begin{equation}
U_i\left( \hat{t}_i \right) =u_i\left( \bar{x}_i,
\hat{t}_i \right) =\phi \left(\frac{\frac{11}{3}(n-1)\hat{t}_i}{b_0}\right)
\left\vert \operatorname{Rm}_{\hat{g}_i} \right\vert ^{2}\left(
\bar{x}_i,\hat{t}_i \right) =1. \label{eq tem use 12}
\end{equation}
On the other hand we have
\begin{align*}
U_i \left( 0\right)  &=\int_{M_i}h_i\left( x,0\right) \left\vert
\operatorname{Rm}_{\hat{g}_i} \right\vert ^{2}\left( x,0\right)
H_i \left( x,0\right) d\mu _{\hat{g}_i\left( 0 \right) } \\
&\leq \int_{B_{\hat{g}_i \left( 0\right) }\left(\bar{ x}_i, 2b_0 \right) }
\left\vert \operatorname{Rm}_{\hat{g}_i}
\right\vert ^{2}\left( x,0\right) H_i \left( x,0\right) d\mu _{\hat{g}_i\left(
0 \right) } \\
&\leq \bar{Q}_i^{-2} \int_{B_{\hat{g}_i\left( 0\right) }\left( \hat{x}_i, 2b_0 \right)
}H_i \left( x,0\right) d\mu _{\hat{g}_i \left( 0\right) } \\
&\leq \bar{Q}_i^{-2} \int_{M_i }H_i \left( x,0\right) d\mu _{\hat{g}_i \left(
0 \right) }
\end{align*}
where we have used support $\operatorname{supp}h_i\left( \cdot ,0\right)
\subset B_{\hat{g}_i\left( 0\right) }\left( \bar{ x}_i, 2b_0 \right)$
in the first inequality and (2iii) in the second inequality.
Hence we have
\begin{equation}
U_i(0) \leq  \bar{Q}_i^{-2}. \label{eq tem use 13}
\end{equation}

By combining (\ref{eq tem use 11}), (\ref{eq tem use 12}), and (\ref{eq tem use 13})
we get
\[
1 \leq e^{\left ( \frac{10}{b_0^2} + 64 \right ) \hat{t}_i}  \bar{Q}_i^{-2}
+\frac{320}{\left ( \frac{10}{b_0^2} + 64 \right )b_0^2}\left(
e^{\left ( \frac{10}{b_0^2} + 64 \right ) \hat{t}_i}-1\right).
\]
Hence
\begin{equation}
1 \leq e^{\left ( \frac{10}{b_0^2} + 64 \right ) (\tilde{\alpha}+1)}  \bar{Q}_i^{-2 }
+\frac{320}{\left ( \frac{10}{b_0^2} + 64 \right )b_0^2} \,
e^{\left ( \frac{10}{b_0^2} + 64 \right )(\tilde{\alpha}+1) }.
\label{eq tem use 14}
\end{equation}
Let
\begin{equation}
b_0 \doteqdot \max \left \{\frac{11}{3} (n-1)(\tilde{\alpha}+1)+1, 3e^{33 
(\tilde{\alpha}+1)} \right \}.  \label{eq tem use 15}
\end{equation}
For such choice of $b_0$ we have 
\[
\frac{320}{\left ( \frac{10}{b_0^2} + 64 \right )b_0^2}
\, e^{\left ( \frac{10}{b_0^2} + 64 \right )(\tilde{\alpha}+1) }< \frac{5}{9}.
\]
Inequality (\ref{eq tem use 14}) is impossible since $\bar{Q}_i \rightarrow \infty$.
Hence we get the required contradiction for Case 2.

\vskip .1cm
\textbf{Case 3} $\tilde{\alpha} = 0$.
The proof for this case is similar to the proof of Case 2.
Let $\hat{t}_i \doteqdot \bar{Q}_i \bar{t}_i$.
Let $\hat{g}_i(t) \doteqdot \bar{Q}_ig_i(\left ( \bar{Q}_i \right )^{-1} t)$,
$t \in [0,\hat{t}_i]$.
By passing to a subsequence we have

\noindent (3i) $\left\vert \operatorname*{Rm}_{\hat{g}_{i}}\right
\vert \left( x,t\right) \leq 4$ for  $x\in B_{\hat{g}_{i}\left( 0\right)
}\left( \bar{x}_{i},A_i\right)$ and $t \in [0,\hat{t}_i]$.

\noindent (3ii) $\left \vert \operatorname*{Rm}_{\hat{g}_{i}} \right \vert
\left(\bar{x}_i, \hat{t}_i \right) =1 $.

\noindent (3iii) $\left\vert \operatorname*{Rm}_{\hat{g}_{i}}\right\vert
\left( x,0\right) \leq \bar{Q}_i^{-1}$
for  $x\in B_{\hat{g}_{i}\left( 0\right) }\left( \bar{x}_{i},A_i\right)$.

\noindent (3iv) $\hat{t}_i \leq \frac{1}{6(n-1)}$, $\hat{t}_i \rightarrow 0$,
$A_i >4e^2$, and $A_i \rightarrow \infty$.

Applying Lemma \ref{lem you choose 2} to $\hat{g}_i(t)$ with $b=2$ we get a function
$h_i: M_i \times [0,\hat{t}_i] \rightarrow [0,1]$ such that the support
\[
\operatorname{supp}h_i(\cdot,t) \subset \overline{B}_{\hat{g}_i(t)} \left
(\bar{x}_i, 4 -\frac{11}{3}(n-1)t \right ) \subset B_{\hat{g}_{i}\left( 0\right)
}\left( \bar{x}_{i},A_i\right)
\]
and
\begin{equation*}
\left( \frac{\partial }{\partial t}-\Delta _{\hat{g}_i\left( t\right) }\right) h_i
\leq \frac{5}{2}h_i.
\end{equation*}

We compute
\begin{align*}
&\left( \frac{\partial }{\partial t}-\Delta _{\hat{g}_i }\right)
\left( h_i\left\vert \operatorname{Rm}_{\hat{g}_i }\right\vert ^{2}\right)  \\
&=\left( \left( \frac{\partial }{\partial t}-\Delta _{\hat{g}_i
}\right) h_i \right) \left\vert \operatorname{Rm}_{\hat{g}_i } \right\vert ^{2}
+h_i\left( \left( \frac{
\partial }{\partial t}-\Delta _{\hat{g}_i }\right) \left\vert \operatorname{Rm
}_{\hat{g}_i }\right\vert ^{2}\right) -2\nabla_{\hat{g}_i } h_i\cdot
\nabla_{\hat{g}_i } \left\vert \operatorname{Rm}_{\hat{g}_i }
\right\vert ^{2} \\
&\leq \frac{5}{2}\, h_i\left\vert \operatorname{Rm}_{\hat{g}_i }\right\vert ^{2}+h_i\left(
-2\left\vert \nabla_{\hat{g}_i } \operatorname{Rm}_{\hat{g}_i }\right\vert^{2}
+16\left\vert \operatorname{Rm}_{\hat{g}_i } \right\vert ^{3}\right)
+ 2 \sqrt{10}\left\vert \operatorname{Rm}_{\hat{g}_i }\right\vert
\cdot h_i^{1/2}\left\vert \nabla_{\hat{g}_i } \operatorname{Rm}_{\hat{g}_i }\right\vert  \\
&\leq \frac{133}{2} \, h_i\left\vert \operatorname{Rm}_{\hat{g}_i }
\right\vert^{2}-2h_i\left\vert \nabla_{\hat{g}_i } \operatorname{Rm}_{\hat{g}_i }
\right\vert ^{2}+ 8 \sqrt{10}\, h_i^{1/2}\left\vert
\nabla_{\hat{g}_i } \operatorname{Rm}_{\hat{g}_i }
\right\vert  \\
&\leq \frac{133}{2} \, \left( h_i\left\vert \operatorname{Rm}_{\hat{g}_i }
\right\vert ^{2}\right) +  80,
\end{align*}
where we have used
\begin{equation*}
\left\vert \nabla_{\hat{g}_i } h_i \right\vert =\frac{\left\vert \phi ^{\prime }\left(
w\right) \right\vert }{2}\left\vert \nabla _{\hat{g}_i }d_{\hat{g}_i\left(
t\right) }\left( x,\bar{x}_i \right) \right\vert _{\hat{g}_i}
\leq \frac{ \sqrt{10}}{2}h_i^{1/2}
\end{equation*}
and $\left\vert \operatorname{Rm}_{\hat{g}_i }\right\vert \leq 4$ on
$\operatorname{supp}h_i(\cdot,t)$.  Here $\phi$ is the function defined in the proof
of Lemma \ref{lem you choose 2}.

Let $u_i \doteqdot h_i\left\vert \operatorname{Rm}_{\hat{g}_i }\right\vert ^{2}$.
We have proved
\[
\left( \frac{\partial }{\partial t}-\Delta _{\hat{g}_i  }\right) u_i
\leq 67 u_i+ 80
\]
on $M_i \times [0,\hat{t}_i]$.

Let $H_i>0$ be the backward heat kernel to the conjugate heat equation
on $(M_i, \hat{g}_i(t) ), \, t \in [0,\hat{t}_i]$,
centered at $\bar{x}_i$. Note that
$\int_{M_i}H_i(\cdot, t) d\mu _{\hat{g}_i\left( t\right) }=1$.
We compute
\begin{align*}
\frac{d}{dt}\int_{M_i}u_iH_id\mu _{\hat{g}_i }
&= \int_{M_i}\left( \left( \frac{\partial }{\partial t}-\Delta _{\hat{g}_i}
\right) u_i \right) H_id\mu _{\hat{g}_i } \\
&\leq \int_{M_i}\left( 67 u_i+ 80 \right) H_id\mu _{\hat{g}_i } \\
&= 67 \int_{M_i}u_i H_i d\mu _{\hat{g}_i }+ 80.
\end{align*}
Hence it follows from a simple integration that $U_i \left( t\right) \doteqdot
\int_{M_i}u_i H_id\mu _{\hat{g}_i}$ satisfies
\begin{equation}
U_i \left( t\right) \leq e^{67t}U_i \left( 0\right) +\frac{80}{67}\left(
e^{67t}-1\right) \label{eq tem use 1}
\end{equation}
for $t \in [0,\hat{t}_i]$.

At $t=\hat{t}_i$ we have
\begin{equation}
U_i\left( \hat{t}_i \right) =u_i\left( \bar{x}_i,
\hat{t}_i \right) =\phi \left(\frac{\frac{11}{3}(n-1)\hat{t}_i}{2}\right)
\left\vert \operatorname{Rm}_{\hat{g}_i} \right\vert ^{2}\left(
\bar{x}_i,\hat{t}_i \right) =1. \label{eq tem use 2}
\end{equation}
On the other hand by an argument similar to the proof of (\ref{eq tem use 13})
we have
\begin{equation}
U_i(0) \leq  \bar{Q}_i^{-2}. \label{eq tem use 3}
\end{equation}

By combining (\ref{eq tem use 1}), (\ref{eq tem use 2}), and (\ref{eq tem use 3})
we get
\[
1 \leq e^{67 \hat{t}_i}  \bar{Q}_i^{-2} +\frac{80}{67}\left(
e^{67 \hat{t}_i}-1\right).
\]
This is impossible since $\hat{t}_i \rightarrow 0$ and $\bar{Q}_i \rightarrow
\infty$. Hence we get the required contradiction for Case 3.

Now we have finished the proof of Proposition \ref{prop variant} modulo the proofs
of Lemma \ref{lem You Are Close} and \ref{lem you choose 2}. 
\qed

\vskip .1cm
\noindent \textbf{3 Proof of two technical lemmas}.
In the proof of Proposition \ref{prop variant} we have used the following
two lemmas. Intuitively the first lemma says that if a ball of radius $1$
has bounded sectional curvature and is volume noncollapsing,
then the isoperimetric constant on small certain size ball
is close to the Euclidean one.
Note that the next lemma and essential the same proof are also given by
Wang \cite{W}.
\begin{lemma} \label{lem You Are Close}
Given $n \geq 2$, $v_{0}>0$ and $\delta_{0}>0$, there is $r>0$, depending
only on $n, \, v_0$, and $\delta_0$,
 which has the following property.
Let $B\left(  x_{0},1\right)$ be a ball in a Riemannian manifold
$\left(M^{n},g \right)$ which satisfies the following:

\noindent (I) The closed $\bar{B}\left(  x_{0},1\right)$ is compact in $M$.

\noindent (II) The Riemann curvature $\left\vert \operatorname{Rm}\right\vert
\leq1$ on $B\left(  x_{0},1\right)$.

\noindent (III) The volume $\operatorname{Vol}\left(  B\left(
x_{0},1\right) \right)  \geq v_{0}>0$.

\noindent Then we have
\begin{equation}
\left ( \operatorname{Area}\left(  \partial\Omega\right) \right ) ^{n}
\geq \left(  1-\delta
_{0}\right)  c_{n}\left(  \operatorname{Vol}\left(  \Omega\right)  \right)
^{n-1} \label{eq isoper first}
\end{equation}
for any regular domain $\Omega\subset B\left(  x_{0},r\right)  $. Here
$c_{n}=n^{n}\omega_{n}$ is the isoperimetric constant for Euclidean space.
\end{lemma}

\noindent {\sl Proof}
{\em Step 1 (Injectivity radius bound)}
Under the assumption of Lemma \ref{lem You Are Close}, by a
theorem of Cheeger-Gromov-Taylor \cite[Theorem A.7]{CGT}
there is a $\iota_{0}>0$ depending only on
$n$ and $v_{0}$ such that the injectivity radius $\operatorname{inj}_{x_0} \geq \iota_{0}$.

{\em Step 2 (Metric tensor on ball $B\left(  x_{0},1 \right)$)}
Let $x=\left(x^{i}\right)$ be the normal coordinates at $x_{0}$.
It follows from a result of Hamilton (see Cao et al \cite[Theorem 4.10, 
page 308]{CCCY})
that for any $\varepsilon>0$ there is $\lambda_{0} = \lambda_{0}\left(
n,\varepsilon\right)$ such that metric tensor
\begin{equation}
\left(  1-\varepsilon\right)  (\delta_{ij})  \leq (g_{ij} ) \leq\left(
1+\varepsilon\right)  ( \delta_{ij} ). \label{eq metric vs euc}
\end{equation}
for $\left\vert x\right\vert \leq \lambda_{0}$.
Note that $\left( \delta_{ij}\right)$ is Euclidean metric in
the coordinates $\left(  x^{i}\right)$.

{\em Step 3 (Approximation argument)}
Let $r\doteqdot\min\left\{  \iota_{0},\lambda_{0}\right\}$
and let $\exp_{x_0}: B(r) \rightarrow B(x_0,r)$ be the exponential map.
$\exp_{x_0}$ is a diffeomorphism. Now we consider
a regular domain $\Omega\subset B\left(  x_{0},r\right)$.
We compute
\begin{align*}
\operatorname{Vol}_{g}\left(  \Omega\right)    & =\int_{\Omega}\sqrt
{\det\left(  g_{ij}\right)  } \cdot dx^{1}\ldots dx^{n}\\
& \leq\int_{\left ( \exp_{x_0} \right )^{-1} \Omega}\sqrt{
\left(  1+\varepsilon\right)  ^{n}\det\left(
\delta_{ij}\right)} \cdot dx^{1}\ldots dx^{n}\\
& =\left(  1+\varepsilon\right)  ^{n/2}\operatorname{Vol}_{\operatorname{Euc}
}\left( \left ( \exp_{x_0} \right )^{-1} \Omega\right)  .
\end{align*}
Let $\{ \theta_{a} \}_{a=1}^{n-1}$ be an orthonormal frame of $\left(
\partial\Omega,\left.  \left(  \delta_{ij}\right)  \right\vert _{\partial
\Omega}\right)  $ at some point $x$ and let $\{ \theta_{a}^{\ast} \}$
be the dual frame. The area form $d\sigma_{\left(  \partial\Omega,\left.  \left(
\delta_{ij}\right)  \right\vert _{\partial\Omega}\right)  }$ at $x$
is given by $\theta_{1}^{\ast}\wedge \ldots
\wedge\theta_{n-1}^{\ast}$. The area form $d\sigma_{\left(  \partial
\Omega,\left. g \right\vert _{\partial\Omega}\right)  }$
at $x$ is given by 
\[
\sqrt{\det\left(  g\left(
\theta_{a},\theta_{b}\right)  \right)  _{\left(  n-1\right)  \times\left(
n-1\right)  }} \cdot \theta_{1}^{\ast}\wedge\ldots\wedge\theta_{n-1}^{\ast}.
\]
We can estimate
\begin{equation*}
\sqrt{\det\left(  g\left(  \theta_{a},\theta_{b}\right)  \right)  _{\left(
n-1\right)  \times\left(  n-1\right)  }}
\geq  \sqrt{\left(  1-\epsilon
\right)  ^{n-1}\det\left( \left(  \delta_{ij}\right)
\left(  \theta_{a},\theta_{b}\right)  \right)}
=\left(  1-\epsilon\right)
^{\left(  n-1\right)  /2},
\end{equation*}
hence
\begin{align*}
& \operatorname{Area}_{\left. g \right\vert _{\partial
\Omega}}\left(  \partial\Omega\right)  =\int_{\partial\Omega}
d\sigma_{\left(  \partial\Omega,\left.  g \right\vert
_{\partial\Omega}\right)  } \\
& \geq \int_{\partial \left ( \left ( \exp_{x_0} \right )^{-1} \Omega \right )
}\left(  1-\varepsilon\right)  ^{\left(
n-1\right)  /2} d\sigma_{\left(  \partial \left (
\left ( \exp_{x_0} \right )^{-1} \Omega \right ),\left.
\left(  \delta_{ij}\right)  \right\vert _{\partial
\left ( \left ( \exp_{x_0} \right )^{-1} \Omega \right ) }\right)  } \\
& =\left(  1-\varepsilon\right)  ^{\left(  n-1\right)  /2}\operatorname{Area}
_{\operatorname{Euc}}\left(  \partial \left (
\left ( \exp_{x_0} \right )^{-1} \Omega \right ) \right)  .
\end{align*}

Now we compute
\begin{align*}
\frac{\left(
\operatorname{Area}_{\left. g \right\vert _{\partial
\Omega}}\left(  \partial\Omega\right)  \right)  ^{n}}{\left(
\operatorname{Vol}_{g}\left(  \Omega\right)  \right)  ^{n-1}}
& \geq\frac{\left(  \left(  1-\varepsilon\right)  ^{\left(  n-1\right)
/2}\operatorname{Area}_{  \operatorname{Euc} }
\left(  \partial \left ( \left ( \exp_{x_0} \right )^{-1} \Omega \right )
\right)  \right)  ^{n}}{\left(
\left(  1+\varepsilon\right)  ^{n/2}\operatorname{Vol}_{\operatorname{Euc}
}\left( \left ( \exp_{x_0} \right )^{-1} \Omega\right)  \right)  ^{n-1}}\\
& =\left(  \frac{1-\varepsilon}{1+\varepsilon}\right)  ^{\frac{\left(
n-1\right)  n}{2}}\cdot\frac{\left(  \operatorname{Area}_{\operatorname{Euc}
}\left(  \partial \left ( \left ( \exp_{x_0} \right )^{-1} \Omega
\right ) \right)  \right)  ^{n}}{\left(  \operatorname{Vol}
_{\operatorname{Euc}}\left( \left ( \exp_{x_0} \right )^{-1}
\Omega\right)  \right)  ^{n-1}}\\
& \geq\left(  \frac{1-\varepsilon}{1+\varepsilon}\right)  ^{\frac{\left(
n-1\right)  n}{2}}c_{n}.
\end{align*}
Given $\delta_0$ we choose
$\varepsilon$ such that
\[
\left(  \frac{1-\varepsilon}{1+\varepsilon}\right)  ^{\frac{\left(
n-1\right)  n}{2}} = 1-\delta_0 ,
\]
this in turn requires us to choose the corresponding
$\lambda_{0}\left(  n,\varepsilon\right)$ to ensure (\ref{eq metric vs euc}).
Then Lemma \ref{lem You Are Close} holds for
$r = \min\left\{  \iota_{0},\lambda_{0}\right\}$.
\qed

\vskip .1cm
The second lemma is about the existence of an auxiliary function.
\begin{lemma} \label{lem you choose 2}
Let $\left( M^{n},g\left( t\right) \right) $, $t\in \left[
0, \hat{t} \right] $, be a solution of Ricci flow. Let
$b$ be a constant bigger than  $\frac{11}{3}(n-1) \hat{t} + 1$
and let $A$  be a constant bigger or equal to $2 e^{4(n-1)\hat{t}} b$.
We assume that closed ball $\overline{B}_{g(0)}(\bar{x},A) \subset M$
be a compact subset and that $| \operatorname{Rm} | (x,t) \leq 4$
for all $(x,t) \in  B_{g(0)} \left (\bar{x}, A \right ) \times [0,\hat{t}]$.
Then there is a function $h: M
\times [0,\hat{t}] \rightarrow [0,1]$ such that
for each $t \in [0,\hat{t}]$ the support
\[
\operatorname{supp}h(\cdot,t) \subset
\overline{B}_{g(t)}(\bar{x}, 2b-\frac{11}{3}(n-1)t) \subset B_{g(0)}(\bar{x},A)
\]
and 
\[
\left( \frac{\partial }{\partial t}-\Delta _{g\left( t\right) }\right) h\leq
\frac{10}{b^{2}} \, h
\]
on $ M \times [0, \hat{t}]$.
\end{lemma}


\noindent {\sl Proof.}
Let $\phi :\mathbb{R}\rightarrow \left[ 0,1\right] $ be a smooth function
which is strictly decreasing on the interval $\left[ 1,2\right] $ and which
satisfies%
\begin{equation}
\phi \left( s\right) =\left\{
\begin{tabular}{ll}
$1$ & if$\text{ }s\in (-\infty ,1],\smallskip $ \\
$0$ & if$\text{ }s\in \lbrack 2,\infty ),$%
\end{tabular}%
\ \ \right.   \label{PhiCutoffInfrared01 add1}
\end{equation}%
and
\begin{subequations}
\label{PhiCutoffInfrared02}
\begin{align}
\left( \phi ^{\prime }(s)\right) ^{2}& \leq 10\phi (s),
\label{eq  temp name 113} \\
\phi ^{\prime \prime }(s)& \geq -10\phi (s)  \label{eq temp name 114}
\end{align}%
for $s\in \mathbb{R}$. We define for any $t\in \left[ 0,T\right] $
\end{subequations}
\begin{equation*}
h\left( x,t\right) =\phi \left( \frac{d_{g\left( t\right) }\left(
x,\bar{x}\right) +at}{b}\right)
\end{equation*}
where $a$ and $b$ are two positive constants to be chosen.
Note that $\operatorname{
supp}h\left( \cdot ,t\right) \subset B_{g\left( t \right) }\left(
\bar{x},2b-at \right) $.

By the curvature assumption we have $B_{g(t)}(\bar{x}, e^{-4(n-1)\hat{t}}A)
\subset B_{g(0)}(\bar{x},A)$ for $t \in [0,\hat{t}]$. We choose
$2b \leq e^{-4(n-1)\hat{t}}A$ so that $\operatorname{
supp}h\left( \cdot ,t\right) \subset B_{g\left( 0 \right) }\left(
\bar{x},A \right) $.

Let $w\left( x,t\right) \doteqdot \frac{d_{g\left(
t\right) }\left( x,\bar{x}\right) +at}{b}$. We compute
\begin{align*}
& \left( \frac{\partial }{\partial t}-\Delta _{g\left( t\right) }\right) h \\
&=\frac{\phi ^{\prime }\left( w\right) }{b}\left( \left( \frac{\partial }{
\partial t}-\Delta _{g\left( t\right) }\right) d_{g\left( t\right) }\left(
x,\bar{x}\right) +a\right) -\frac{\phi^{\prime\prime}\left( w\right) }{b^{2}}\left\vert
\nabla _{g\left( t\right) }d_{g\left( t\right) }\left( x,\bar{x}\right)
\right\vert _{g\left( t\right) }^{2} \\
& \leq \frac{\phi ^{\prime }\left( w\right) }{b}\left( \left( \frac{\partial }{
\partial t}-\Delta _{g\left( t\right) }\right) d_{g\left( t\right) }\left(
x,\bar{x}\right) +a \right ) +\frac{10}{b^{2}} \, h.
\end{align*}
Choosing $a$ such that $a \hat{t} < b -1$, then for $x \in B_{g(t)}(\bar{x},1)$
or $x \notin \operatorname{supp} h(\cdot,t)$  we have
$\phi ^{\prime }\left( w\right)(x,t) =0$. Hence for such $x$ we have
\[
\left( \frac{\partial }{\partial t}-\Delta _{g\left( t\right) }\right) h\leq
\frac{10}{b^{2}} \, h.
\]

For $x \notin B_{g(t)} \left(\bar{x}, 1 \right)$ and $x \in \operatorname{supp}
h(\cdot,t)$,
we use  \cite[Lemma 8.3(a)]{Pe02I} with $r_{0}=1$ and $K=4$ and get
\begin{equation*}
\left. \left( \frac{\partial }{\partial t}-\Delta _{g(t)}\right)
d_{g(t)}(x,\bar{x})\right\vert _{t=t_{0}}\geq -(n-1)\left( \frac{2}{3}Kr_{0}+
\frac{1}{r_{0}}\right) =-\frac{11}{3}(n-1).
\end{equation*}
By choosing $a \doteqdot \frac{11}{3}(n-1)$ and using $\phi^\prime (w) \leq 0$
we obtain
\begin{equation*}
\left( \frac{\partial }{\partial t}-\Delta _{g\left( t\right) }\right) h\leq
\frac{10}{b^{2}} \, h.
\end{equation*}
The lemma is proved.
\qed

\vskip .1cm
\noindent \textbf{4 Two examples}.  In this section we give two example 
showing that neither
the volume lower bound assumption
nor the completeness assumption in Theorem \ref{thm variant} can be dropped.

Let $r$ be an arbitrary positive constant in $(0,1]$.
Let $g_r^0$ be a Riemannian metric on a
topological sphere $\Sigma^2$ which contains a round cylinder
$S^1(r) \times [-1,1]$ of radius $r$ and length $2$.
We have $\operatorname{Vol}_{g^0_r}(\Sigma) \geq 4 \pi r$.
We assume volume $\operatorname{Vol}_{g^0_r}(\Sigma) \leq 20r$.
Let  $(\Sigma^2, g_r(t))$, $t \in [0,T_r)$, be the maximal solution of
the Ricci flow with $g_r(0)=g_r^0$. Then the blowup time
\[
T_r = \frac{1}{8 \pi} \operatorname{Vol}_{g^0_r}(\Sigma)
\in  \left ( \frac{1}{2}r,  \frac{5}{2 \pi}r \right ].
\]

Let $p \in S^1(r)$. Then $x_0 \doteqdot (p,0)$ is a point in $\Sigma$.
For any $\epsilon_0$ we can choose $r$ small enough so that $T_r < \epsilon_0$.
Clearly we have $\left \vert \operatorname{Rm}_{g_r} \right \vert (x,0) =0$ for
 $x \in B_{g_r(0)}(x_0,1)$ and $\operatorname{Vol}_{g_r(0)}(B_{g_r(0)}(x_0,1))
\leq 4\pi r$. For any $\epsilon \in (\frac{1}{2}r,T_r)$, should the conclusion
of Theorem \ref{thm variant} hold for $g_r(t)$ when $r$ is small enough,
we would have
$\left \vert \operatorname{Rm}_{g_r} \right \vert (x_0,\epsilon) \leq
\epsilon_0^{-2}$. Since $\epsilon$ is arbitrary, we have $\lim_{t \rightarrow T_r}
\left \vert \operatorname{Rm}_{g_r} \right \vert (x_0,t) < \epsilon_0^{-2}$. However
it is well-known that the limit should be infinity. Hence Theorem \ref{thm variant}
does not hold for $g_r(t)$.
By taking the product of $(\Sigma^2, g_r(t))$ with flat torus we get high dimensional
examples.

The second example is a simple modification of the previous example,
the idea of construction is due to Peter Topping (unpublished work).
Let
\[
\Phi: \mathbb{R} \times (-1,1) \rightarrow S^1(r) \times (-1,1) \subset \Sigma
\]
be the standard universal cover map. Then $\left( \mathbb{R} \times (-1,1),
\Phi^*g_r(t) \right )$ is a incomplete solution of the Ricci flow.
Clearly we have $\left \vert \operatorname{Rm}_{\Phi^*g_r} \right \vert (x,0) =0$
for $x \in B_{\Phi^*g_r(0)}((0,0),1)$ and $\operatorname{Vol}_{\Phi^*g_r(0)}(B_{\Phi^*g_r(0)}((0,0),1))= \pi $.
Arguing as in the previous example we conclude that
Theorem \ref{thm variant} does not hold for $\Phi^*g_r(t)$ with the ball
center being $(0,0)$  when $r$ is small enough.

Finally we make a remark. It follows from the proof of Lemma \ref{lem You Are Close}
that under the same assumption as the lemma there is a $\tilde{r}\in (0,1]$,
depending only on $n$, $v_0$, and $\delta_0$, such that
\[
\operatorname{Vol} \left ( B(x_0, \tilde{r} ) \right )
\geq (1-\delta_0 ) \omega_n \tilde{r}^n.
\]
We need to switch the notations below. Denote the $r_0$ in
Theorem \ref{Thm 10.3} by $r_1$ and denote the $r_0$ in Theorem \ref{thm variant}
still by $r_0$. Let $\delta_0$ be the $\delta$ in Theorem \ref{Thm 10.3}.
Let $g(t)$ be a solution of the Ricci flow satisfying the assumption of Theorem
\ref{thm variant}. Then the assumption of Theorem \ref{Thm 10.3} holds
for $g(t)$ with $r_1 = r_0 \tilde{r}$, hence by Theorem \ref{Thm 10.3}
we get a curvature bound which is essentially equivalent to the curvature
bound given by Theorem \ref{thm variant}.
The reason, why we do not use Theorem \ref{Thm 10.3} and the proof of
Lemma \ref{lem You Are Close} to give a more direct proof of
Theorem \ref{thm variant} is that at the time of writing this note
the author is not aware of a detailed proof
of Theorem \ref{Thm 10.3} in the literature.

\section*{Acknowledgments}
{Part of this work was done while the author was visiting Beijing Center for
Mathematical Research during the turn of 2008 to 2009 and Department of Mathematics
at University of California at San Diego in early 2009.
The author thanks Professor Weiyue Ding, Gang Tian, and Bennett Chow for
their invitation and hospitality. The author also thanks the referee for suggestions
in improving the presentation of this note.}

\bibliographystyle{natbib}






\medskip
\noindent {\sc Peng Lu}, University of Oregon

\noindent e-mail: penglu@uoegon.edu

\end{document}